# The Riccati Characteristic Equation

Douglas R. Frey, *Fellow, IEEE*

**Abstract:** The Riccati differential equation is examined in light of its connection to second order linear time varying systems. In that light it becomes the clear generalization for the characteristic equation of linear time invariant systems, and is called the Riccati Characteristic Equation (RCE). Consequently, the RCE becomes the unifying centerpiece for the study of linear systems. Its solutions are considered in complementary pairs that form a continuum based on a primitive pair. Pairs may always be found as purely real solutions, despite the fact that complex conjugate primitive solutions are shown to exist in many cases. Not only is the pairing unique, but the general form of solutions, shown here for the first time, is uniquely compact and encompasses all known solutions, while allowing for all initial conditions. Classical engineering mathematics examples are shown to conform to this approach, which provides new insights to all, especially Floquet theory.

Index terms–Control systems, Eigenvalues and eigenfunctions, Linear systems, Mathematical analysis, Nonlinear equations, Riccati equations, State-space methods, Time-varying systems[1].

## I. Introduction

The Riccati equation in its various forms has found use in many areas of engineering, especially control systems. From pure mathematics to the analysis of time varying systems, this equation has proven to be of interest. Unfortunately, the solution to the scalar Riccati differential equation covered here has been daunting. Its connection to linear time varying dynamical systems, while well known, has not been fully exploited. This paper and its companion paper [1] addresses that by giving the Riccati differential equation a central role in the behavior of linear systems. By offering the Riccati equation as the Riccati characteristic equation, it generalizes the well known characteristic equation in a way suitable for all time varying systems, without interfering with its usage in time invariant systems. Consequently, this paradigm of analysis may impact the analysis and design of time varying systems across the full range of applications, from adaptive filters to communications electronics, and in many other areas of engineering as well. Moreover, its usage in the study of engineering mathematics such as the Bessel, quantum harmonic oscillator, and Matthieu equations, as shown in this work, provides potentially powerful new insight.

Specifically, this paper investigates what has classically been called the Riccati Differential Equation, which in the present context is now referred to as the Riccati Characteristic Equation (RCE) due to its intimate connection to linear time varying (LTV) systems and their solution, particularly as described in the companion paper [1]. We assume that all Riccati equation coefficients are real, as would be expected for the associated second order practical systems. This paper offers a new paradigm regarding many aspects of the Riccati equation. For example, the idea that finite escape time precludes the existence of a solution has often been suggested [e.g.2]. The current work not only admits finite escape time behavior, but shows that such behavior is a very natural part of solutions, owing to the well known RCE connection to linear systems.

The explicit nature of this connection and the solutions of the RCE are at the heart of this work. Section II reviews the mathematical basis and introduces appropriate preliminary material. Section III makes the logical connection of the Riccati equation to the general characteristic equation for a LTV second order system–that is, as the Riccati Characteristic Equation, or RCE. A natural consequence of this connection is that solutions to the RCE are studied in pairs, with primitive pairs helping to define the general forms of solutions to the RCE. Such simple general solution forms for the Riccati equation have never been given before. Section IV discusses the utility of these pairs of solutions in finding general solutions to the associated second order LTV systems, providing a completely intuitive extension of classic LTI system solutions to LTV systems. The work in [1] fully articulates this. Section V shows how the paradigm introduced here regarding solutions offers new insight into the complete unification of known solutions for the Riccati equation, the phase portrait of the Riccati equation, and the fundamental existence of solutions to the Riccati equation, including periodic solutions. In particular, the work here and in [1] show for the first time that periodic solutions must exist for Riccati equations with periodic coefficients. Finally, solutions to well studied systems, specifically the Bessel, quantum harmonic oscillator, and Matthieu equations are also shown to be just special cases of the solutions shown here. New perspectives on the solutions of these systems are highlighted in each case.

## II. Mathematical Preliminaries

We begin by reviewing the mathematical basis. The equation below represents the most general form of the Riccati equation considered here.

$$\dot{z}(t) = s_2(t)z^2(t) + s_1(t)z(t) + s_0(t) \qquad (1)$$

The solution of this equation has proven to be elusive. While

---

The author is with the ECE Dept., Lehigh University, Bethlehem, PA 18015 (e-mail: drf3@lehigh.edu)

many techniques have been tried, we will begin by focusing on an approach that has often been used to simplify things [3,7]. Specifically, one may use a change of variables to simplify the equation to the form below for v(t), as shown in (2), where η(t) = -$s_1$(t)/(2$s_2$(t)).

$$\dot{z}(t) = s_2(t)z^2(t) + s_1(t)z(t) + s_0(t)$$
$$= s_2(t)(z(t) + \frac{s_1(t)}{2s_2(t)})^2 - s_2(t)(\frac{s_1(t)}{2s_2(t)})^2 + s_0(t)$$
$$\dot{z}(t) - \dot{\eta}(t) = -\omega_{01}(t)(z(t) - \eta(t))^2 + \omega_{02}(t) \quad (2)$$
$$\omega_{01}(t) = -s_2(t) \; ; \; \omega_{02}(t) = s_0(t) - s_2(t)(\frac{s_1(t)}{2s_2(t)})^2 - \dot{\eta}(t)$$
$$\dot{\nu}(t) = -\omega_{01}(t)\nu^2(t) + \omega_{02}(t) \; ; \; \nu(t) = z(t) - \eta(t)$$

The final reduced form of the Riccati equation for v(t) has been chosen to match that in the discussion of LTV system analysis in [1]. This specifically relates to the first step of the state space system transformation process using $M_0$(t) in [1].

The scalar related LTV system context here is that of a second order differential equation. Specifically, consider the equation below.

$$\ddot{y}(t) + r_1(t)\dot{y}(t) + r_0(t)y(t) = 0 \quad (3)$$

It is of interest to note that a companion form for this differential equation produces these state equations.

$$x_1(t) = y(t) \; ; \; x_2(t) = \dot{x}_1(t) = \dot{y}(t)$$
$$\begin{bmatrix} \dot{x}_1(t) \\ \dot{x}_2(t) \end{bmatrix} = \begin{bmatrix} 0 & 1 \\ -r_0(t) & -r_1(t) \end{bmatrix} \begin{bmatrix} x_1(t) \\ x_2(t) \end{bmatrix} \quad (4)$$

The state matrix in (4) is a special case of the most general LTV state matrix, A(t), for a second order system. Using the derivation in [1], this special case, and the general case, possess the Riccati characteristic equation for v(t) shown in (5), where $a_{ij}$ is the time-varying ij-th element of A(t). The specifics relating to the special case in (4) are also given, where all parameters are time varying.

$$\dot{\nu} = -\omega_{01}\nu^2 + \omega_{02} \; ; \; \nu = z - \eta$$
$$\omega_{01} = a_{12} = 1 \; ; \; \eta = \frac{a_{11} + a_{22}}{2a_{12}} = -\frac{1}{2}r_1$$
$$\alpha = \frac{a_{11} - a_{22}}{2a_{12}} = \frac{1}{2}r_1 \quad (5)$$
$$\omega_{02} = \dot{\alpha} + a_{12}\alpha^2 + a_{21} = \frac{1}{2}\dot{r}_1 + \frac{1}{4}r_1^2 - r_0$$

We are now in a position to analyze the Riccati equation in the reduced form shown without loss of generality. Because of the duality of a linear second order system and its associated Riccati characteristic equation, given by how the parameters, $\omega_{01}$(t), $\omega_{02}$(t), and η(t), are related to the components of the state matrix, we will have the need to consider RCE solutions as intrinsically tied to those of the associated LTV system. To that end, we will refer to solutions of the LTV system as being in "the time domain" and solutions to the associated Riccati characteristic equation as being solutions in "the Riccati domain". Note that both types of solutions will be written using time, t, as the independent variable, since the duality ensures that they track at each instant in time. This is highlighted by the fact that each given solution, y(t), in the time domain is related to its counterpart, z(t), in the Riccati domain as shown below. The dynamic eigenvalue, λ(t), which is a scaled version of z(t), explicitly defines the time domain solution.

$$y(t) = time\ domain\ solution = Ae^{\int \lambda(t)dt}$$
$$z(t) = Riccati\ domain\ solution = \nu(t) + \eta(t) \quad (6)$$
$$\lambda(t) = \omega_{01}(t)z(t)$$

It is noted that the time domain solution in this discussion is, in general, not the explicit solution desired in a given case, since the given Riccati domain solution directly specifies one of the state variables, $x_1$(t) in the present case, in the state space description of the second order system. Nevertheless, we may assume that the core components of the time domain solution we use here are those of any explicit solution desired in practice except in trivial special cases, where for example the state equations are uncoupled. In fact, every possible real time domain solution for $x_1$(t) to within a scale factor, with specific boundary conditions, is explicitly given by one of the solutions to the associated RCE, which will be shown later. In addition, the broad range of solutions to the LTV system in the time domain are completely specified using the fundamental matrix for the system, which is also completely specified using the pair of solutions of the RCE.

**III. The Riccati Characteristic Equation and Solutions**

The use of the Riccati differential equation as the characteristic equation for all second order linear systems is motivated by the classic linear time invariant (LTI) case, and has been proposed before [4-6], but not in the way here. The Riccati differential equation becomes the usual quadratic algebraic equation, with a pair of solutions, when the derivative is assumed to be zero, as in the LTI case. However, in the LTV case, the solutions to the Riccati equation may be examined in pairs as well. As shown below, the Riccati characteristic equation always has primitive solutions that may be put in the form, $v_R$(t) ± $v_I$(t), where $v_R$(t) is real and $v_I$(t) is either real or purely imaginary. This offers a uniformity to the analysis across the full spectrum of LTV systems. The primitive solutions also provide important aspects of the phase portrait of the Riccati equation.

Classic LTI analysis confirms that one may always solve the characteristic equation by inserting the guess, $\lambda_R \pm \lambda_I$, into the characteristic equation, where $\lambda_I$ may be real or imaginary. Specifically, we have,

$$\begin{aligned} 0 &= a\lambda^2 + b\lambda + c \\ &= a(\lambda_R \pm \lambda_I)^2 + b(\lambda_R \pm \lambda_I) + c \\ &= a\lambda_R^2 + a\lambda_I^2 + b\lambda_R + c \pm (2a\lambda_R\lambda_I + b\lambda_I) \quad (7) \\ \Rightarrow 0 &= a\lambda_R^2 + a\lambda_I^2 + b\lambda_R + c \; ; \; 0 = 2a\lambda_R\lambda_I + b\lambda_I \\ \lambda_R &= -\frac{b}{2a} \; ; \; \lambda_I = \sqrt{\frac{b^2}{4a^2} - \frac{c}{a}} \end{aligned}$$

The final result, when combined to produce $\lambda = \lambda_R \pm \lambda_I$, is no surprise. These steps, however, suggest the same analysis in the LTV case, where we again assume the same form for the pair of solutions to the Riccati equation in (1) and (2), directly related to the dynamic eigenvalues. All quantities in (8) are now time varying.

$$\begin{aligned} \dot{z}_R \pm \dot{z}_I &= s_2(z_R \pm z_I)^2 + s_1(z_R \pm z_I) + s_0 \\ &= s_2 z_R^2 + s_2 z_I^2 + s_1 z_R + s_0 \pm (2s_2 z_R z_I + s_1 z_I) \\ \dot{z}_R &= s_2 z_R^2 + s_2 z_I^2 + s_1 z_R + s_0 \; ; \; \dot{z}_I = 2s_2 z_R z_I + s_1 z_I \quad (8) \\ z_R &= -\frac{s_1}{2s_2} + \frac{1}{2s_2}\frac{\dot{z}_I}{z_I} \end{aligned}$$

As in the analysis in (2), we may reduce and simplify the form of this last result, and in the process allow uniformity of the discussion. This exploits the fact that z(t) and v(t) are just shifted versions of one another.

$$\begin{aligned} \dot{\nu} &= -\omega_{01}\nu^2 + \omega_{02} \\ \nu &= \nu_R \pm \nu_I = z_R \pm z_I - \eta \; ; \; \eta = -\frac{s_1}{2s_2} \\ \omega_{01} &= -s_2 \; ; \; \omega_{02} = s_0 - s_2\eta^2 - \dot{\eta} \quad (9) \\ \nu_R &= z_R - \eta = \frac{1}{2s_2}\frac{\dot{z}_I}{z_I} \; ; \; \nu_I = z_I \\ \therefore \nu_R &= -\frac{1}{2\omega_{01}}\frac{\dot{\nu}_I}{\nu_I} \Rightarrow \nu = -\frac{1}{2\omega_{01}}\frac{\dot{\nu}_I}{\nu_I} \pm \nu_I \end{aligned}$$

It is noteworthy that what we shall call the *real component*, $\nu_R(t)$, of the RCE solution is a function of what we shall call the *intrinsic component*, $\nu_I(t)$. This form of the pair of solutions is noteworthy, but the existence of the pair is not. It has been long understood that if one solution of the Riccati equation were known, then a second solution could always be found. In the present context, given that v(t) is a solution of the equation, then v(t) + κ(t) must also be a solution for the correct function κ(t). This may be proven by substitution into the Riccati equation and solving the Bernoulli equation that results [e.g. 7]. A slight variation on this idea may be had by assuming that if v(t) = v_R(t) + v_I(t) is a solution of the equation, then v(t) - 2/p(t) must also be a solution for the correct function, p(t), whose utility is shown in [1]. Given that v(t) solves the equation, let us substitute v(t) - 2/p(t) into the RCE to see the constraint on p(t).

$$\frac{d}{dt}(\nu(t) - \frac{2}{p(t)}) = -\omega_{01}(t)(\nu(t) - \frac{2}{p(t)})^2 + \omega_{02}(t)$$
(10)
$$\Rightarrow \dot{p}(t) = 2\omega_{01}(t)\nu(t)p(t) - 2\omega_{01}(t)$$

Clearly, this is a linear first order differential constraint which allows an analytic solution for p(t). It is valuable that this is the equation for p(t) that enables the spectrum invariant transformation of the state matrix to be established in [1]. The particular solution to this equation is given below.

$$\begin{aligned} p(t) &= e^{\int 2\omega_{01}(t)\nu(t)dt} \int (-2\omega_{01}(t))e^{-\int 2\omega_{01}(t)\nu(t)dt} dt \\ \nu(t) &= \nu_R(t) + \nu_I(t) = -\frac{\dot{\nu}_I(t)}{2\omega_{01}(t)\nu_I(t)} + \nu_I(t) \\ \int 2\omega_{01}(t)\nu(t)dt &= -\ln(\nu_I(t)) + \int 2\omega_{01}(t)\nu_I(t)dt \\ e^{\int 2\omega_{01}(t)\nu(t)dt} &= \frac{1}{\nu_I(t)}e^{\int 2\omega_f(t)dt} \; ; \; \omega_f(t) = \omega_{01}(t)\nu_I(t) \quad (11) \\ p(t) &= \frac{e^{\int 2\omega_f(t)dt}}{\nu_I(t)} \int (-2\omega_f(t))e^{-\int 2\omega_f(t)dt} dt \\ &= \frac{e^{\int 2\omega_f(t)dt}}{\nu_I(t)} e^{-\int 2\omega_f(t)dt} = \frac{1}{\nu_I(t)} \end{aligned}$$

Not only does this specify a formula for p(t), but it completely justifies the complementary solution. Specifically, given that p(t) = 1/v_I(t), then -2/p(t) = -2v_I(t), and therefore v(t) - 2/p(t) = v(t) - 2v_I(t). Given that v(t) = v_R(t) + v_I(t) is a solution to the Riccati equation, then v_R(t) + v_I(t) - 2v_I(t) = v_R(t) - v_I(t) must also be a solution, proving the consistency of all the above computations.

Under certain circumstances, the pair of solutions just considered may be called "primitive" solutions. However, there is a continuum of solutions to the Riccati equation, since one need not use the particular solution for p(t) in obtaining the "second" solution. Specifically, we may use the complete solution for p(t), labeled p_c(t) below, which introduces a constant variable of integration, called $K_0$ below. Following the calculations above, we have,

$$p_c(t) = \frac{e^{2\varphi_f(t)}}{\nu_I(t)}[e^{-2\varphi_f(t)} + K_0] \; ; \; \varphi_f(t) = \int \omega_f(t)dt$$
$$\nu_{Ic}(t) = \frac{1}{p_c(t)} = \nu_I(t)\frac{e^{-2\varphi_f(t)}}{e^{-2\varphi_f(t)} + K_0} \quad (12)$$

We may now find a new second solution, $\nu_2(t)$, for v(t) by subtracting $2/p_c(t)$ from the first solution, called $\nu_1(t)$ below.

$$\nu_2(t) = \nu_1(t) - 2\nu_{1c}(t) = \nu_R(t) + \nu_I(t) - 2\nu_{1c}(t)$$
$$= -\frac{\dot{\nu}_I(t)}{2\omega_{01}(t)\nu_I(t)} + \nu_I(t) - 2\nu_I(t)\frac{e^{-2\varphi_f(t)}}{e^{-2\varphi_f(t)} + K_0} \quad (13)$$
$$= -\frac{\dot{\nu}_I(t)}{2\omega_{01}(t)\nu_I(t)} - \nu_I(t)\frac{e^{-2\varphi_f(t)} - K_0}{e^{-2\varphi_f(t)} + K_0}$$

Naturally, setting $K_0 = 0$ recovers the complementary solution given earlier. An interesting pair of possibilities result from choosing $K_0$ equal to plus and minus 1.

$$\nu_2(t) = -\frac{\dot{\nu}(t)}{2\omega_{01}(t)\nu_I(t)} - \nu_I(t)\frac{e^{-2\varphi_f(t)} - 1}{e^{-2\varphi_f(t)} + 1} \; ; \; K_0 = 1$$
$$= -\frac{\dot{\nu}(t)}{2\omega_{01}(t)\nu_I(t)} + \nu_I(t)\tanh(\varphi_f(t)) \quad (14)$$
$$\nu_2(t) = -\frac{\dot{\nu}(t)}{2\omega_{01}(t)\nu_I(t)} + \nu_I(t)\coth(\varphi_f(t)) \; ; \; K_0 = -1$$

Other values of $K_0$ lead to a simple generalization of these results. In order to see this, note the following, covering cases where $K_0$ is implicitly positive or negative and replaced by its absolute value. K is an assumed constant of integration in finding $\varphi_f(t)$, such that $K = -\frac{1}{2} \ln(|K_0|)$.

$$\frac{e^{-2\varphi_f(t)} - K_0}{e^{-2\varphi_f(t)} + K_0} = \frac{K_0 e^{-2(\varphi_f(t)-K)} - K_0}{K_0 e^{-2(\varphi_f(t)-K)} + K_0}$$
$$= -\tanh(\varphi_f(t) - K) \quad (15)$$
$$\frac{e^{-2\varphi_f(t)} + K_0}{e^{-2\varphi_f(t)} - K_0} = -\coth(\varphi_f(t) - K)$$

Now observe that if the tanh solution of (14) and (15) is used as a first known solution, it is easily shown that the complementary second solution, using the corresponding particular solution for p(t), is the coth solution. Furthermore, in each solution, an infinite value of K defaults the solution to the primitive one without the hyperbolic functions. Thus we have the following pair of complementary solutions as representative of the continuum of solutions to the Riccati characteristic equation.

$$\nu_1(t) = -\frac{\dot{\nu}_I(t)}{2\omega_{01}(t)\nu_I(t)} + \nu_I(t)\tanh(\varphi_f(t) - K)$$
$$\nu_2(t) = -\frac{\dot{\nu}_I(t)}{2\omega_{01}(t)\nu_I(t)} + \nu_I(t)\coth(\varphi_f(t) - K) \quad (16)$$

For completeness, it is worth noting that the continuum of solutions creates a continuum of complementary pairs of solutions. Specifically, the pair in (16) obeys all the properties of the primitive pair of solutions. For example, we may write the following.

$$\text{Let } \nu_1(t) = \nu_{Rx}(t) + \nu_{Ix}(t) \; , \; \nu_2(t) = \nu_{Rx}(t) - \nu_{Ix}(t)$$
$$\nu_{Ix}(t) = \frac{\nu_I(t)}{2}(\tanh(\varphi_f(t) - K) - \coth(\varphi_f(t) - K))$$
$$\nu_{Rx}(t) = -\frac{\dot{\nu}_{Ix}(t)}{2\omega_{01}(t)\nu_{Ix}(t)} = -\frac{\dot{\nu}_I(t)}{2\omega_{01}(t)\nu_I(t)} + \quad (17)$$
$$\frac{\nu_I(t)}{2}(\tanh(\varphi_f(t) - K) + \coth(\varphi_f(t) - K))$$

The final result verifies that the general real component, $\nu_{Rx}(t)$, is simply ½ times the sum of the 2 complementary solutions in (16), which is similarly related, as in (9), to what we may call the general intrinsic component, $\nu_{Ix}(t)$, as given in (17). Adding $\pm\nu_{Ix}(t)$ to $\nu_{Rx}(t)$ yields both members of the complementary pair as shown in (16). Furthermore, direct substitution of the solutions in (16) into the reduced RCE proves that all pairs shown solve the RCE, as long as the primitive pair, $v(t) = \nu_R(t) \pm \nu_I(t)$, solves the equation.

It should be noted that the above computations implicitly assume that the primitive pair of solutions, $\nu_R(t) \pm \nu_I(t)$, is purely real, since a real value of $K_0$ is required to ensure that the time domain solutions associated with the Riccati domain solutions can be real via any appropriate choice of constant scaling. We choose to ignore any Riccati domain solution that cannot result in a purely real time domain solution. Nevertheless, another variation on the solutions of (16) is due to the possibility that the primitive solutions $\nu_R(t) \pm \nu_I(t)$ are complex conjugates, which can happen if $\nu_I(t)$ is purely imaginary, since $\nu_R(t)$ still remains real due to its relation to $\nu_I(t)$. A complex value for $\nu_I(t)$ precludes real associated time domain solutions, which will again be ignored. In the case of an imaginary value for $\nu_I(t)$, we may write,

$$\nu_I(t) = j\nu_{\text{Im}}(t) \;;\; \omega_f(t) = j\omega_{fm}(t) \,,\, \varphi_f(t) = j\varphi_{fm}(t)$$
$$\nu_I(t)\tanh(\varphi_f(t) - K) = -\nu_{\text{Im}}(t)\tan(\varphi_{fm}(t) - K)$$
$$\nu_I(t)\coth(\varphi_f(t) - K) = \nu_{\text{Im}}(t)\cot(\varphi_{fm}(t) - K) \quad (18)$$
$$\nu_1(t) = -\frac{\dot{\nu}_{\text{Im}}(t)}{2\omega_{01}(t)\nu_{\text{Im}}(t)} - \nu_{\text{Im}}(t)\tan(\varphi_{fm}(t) - K)$$
$$\nu_2(t) = -\frac{\dot{\nu}_{\text{Im}}(t)}{2\omega_{01}(t)\nu_{\text{Im}}(t)} + \nu_{\text{Im}}(t)\cot(\varphi_{fm}(t) - K)$$

The constant of integration, K, has been associated with the real integral in the above calculations. Thus, we now have compact pairs of purely real solutions to the RCE for the continuum of cases given either real or complex conjugate pairs of primitive solutions. This yields satisfying real solutions to the associated systems in the time domain. The primitive solutions can yield complex conjugate solutions in the time domain, but scaling with complex conjugate constants will produce real solutions as is well known.

The above calculations in (18) show that the solutions in (16) actually still apply for the complex conjugate case if we simply allow a purely imaginary value for K. Therefore, for the first time the general form for solutions of the differential Riccati equation may be compactly written. Specifically, using the results in (2) and (16), we may now write that any practical solution (related to the RCE), z(t), to the Differential Riccati equation, repeated for convenience below, may be written as shown next.

$$\dot{z}(t) = s_2(t)z^2(t) + s_1(t)z(t) + s_0(t)$$
$$z(t) = \nu(t) + \eta(t) \;;\; \eta(t) = -\frac{s_1(t)}{2s_2(t)}$$
$$= \eta(t) - \frac{\dot{\nu}_I(t)}{2\omega_{01}(t)\nu_I(t)} + \nu_I(t)\begin{cases}\tanh(\varphi_f(t) - K) \\ \coth(\varphi_f(t) - K)\end{cases} \quad (19)$$
$$\omega_{01}(t) = -s_2(t) \;;\; \varphi_f(t) = \int \omega_{01}(t)\nu_I(t)dt$$

This solution form is completely contingent upon the intrinsic component of the solution to the reduced equation we have called the RCE, whose primitive pair of solutions is given by $\nu_R(t) \pm \nu_I(t)$.

Clearly, if the primitive pair is known then the general solution is completely specified. However, only one solution of the RCE may be known, that may or may not be primitive. Nevertheless, we may decompose any known solution into its real and intrinsic parts. Suppose the known solution, v(t), is as given in (16) with the tanh option, where $K = -\frac{1}{2}\ln(K_0)$ from earlier. The complete solution, $p_c(t)$, may now be found, where $K_1$ is the new constant of integration

$$\nu(t) = -\frac{\dot{\nu}_I(t)}{2\omega_{01}(t)\nu_I(t)} + \nu_I(t)\tanh(\varphi_f(t) - K)$$
$$\Rightarrow e^{\int 2\omega_{01}(t)\nu(t)dt} = \frac{1}{\nu_I(t)}\cosh^2(\varphi_f(t) - K)$$
$$p_c(t) = \frac{1}{\nu_I(t)}\cosh^2(\varphi_f(t) - K)\int \frac{-2\omega_{01}(t)\nu_I(t)}{\cosh^2(\varphi_f(t) - K)}dt \quad (20)$$
$$= \frac{-2}{\nu_I(t)}\cosh^2(\varphi_f(t) - K)(\tanh(\varphi_f(t) - K) - \frac{K_1}{2})$$
$$= \frac{-2}{\nu_I(t)}\cosh(\varphi_f(t) - K)\sinh(\varphi_f(t) - K)w(K_1, t)$$
$$w(K_1, t) = 1 - \frac{K_1}{2}\coth(\varphi_f(t) - K)$$

The particular solution, where $K_1 = 0$, that will take tanh to coth in the RCE solution by subtracting $2/p_c(t)$, is the same as before. A careful exercise in algebra and trigonometry shows that the full continuum of solutions may be obtained with appropriate values of $K_1$. It is noteworthy that by setting $K_1$ to -2, the primitive solution will be obtained. $p_c(t)$ for that case is shown in the first line below.

$$p_c(t) = \frac{-2w(-2, t)}{\nu_I(t)}\cosh(\varphi_f(t) - K)\sinh(\varphi_f(t) - K)$$
$$-\frac{2}{p_c(t)} = \frac{\nu_I(t)(\coth(\varphi_f(t) - K) - \tanh(\varphi_f(t) - K))}{1 + \coth(\varphi_f(t) - K)} \quad (21)$$
$$= \nu_I(t)(1 - \tanh(\varphi_f(t) - K))$$

Clearly, if $2/p_c(t)$ is subtracted from v(t) we obtain the primitive solution. Had v(t) included the coth option, the same computations would result in a primitive solution. Finally, given a primitive solution, it has already been shown that finding the particular solution for p(t) yields the intrinsic solution as $\nu_I(t) = 1/p(t)$. Thus, we have the following method of finding the intrinsic solution, and therefore, the general solution of any system, given just one arbitrary solution. We assume to begin that the primitive solution is real, in the form of (16), to aid understanding.

$$g(t) = \int (-2\omega_{01}(t))e^{-\int 2\omega_{01}(t)\nu(t)dt}dt$$
$$\nu_p(t) = \nu(t) + \frac{1}{\omega_{01}(t)}\frac{\dot{g}(t)}{g(t) - 2} = \nu_R(t) + \nu_I(t)$$
$$g_p(t) = \int (-2\omega_{01}(t))e^{-\int 2\omega_{01}(t)\nu_p(t)dt}dt \quad (22)$$
$$\nu_I(t) = \frac{-\dot{g}_p(t)}{2\omega_{01}(t)g_p(t)}$$

Note that if the given solution is real, but corresponds to a complex primitive solution, the integral for g(t) could technically yield an imaginary function, if the imaginary intrinsic solution, $v_I(t) = jv_{Im}(t)$, were taken into account in evaluating the exponential of the integral of $v_R(t)$. With that, g(t) and it's derivative would be imaginary and the choice of $K_1$ equals -2 yields the correct complex primitive solution. This requires the correct guess of the nature of the primitive solution in using (22). However, the nature of the given solution will likely make this guess straightforward, especially in light of the discussion below. This now completely specifies the RCE, its solutions, and the connection between solutions.

## IV. Riccati and Time Domains Solutions

Using these Riccati equation results, the general solution for the dynamic eigenvalues for the associated LTV system, as expressed in (6) where $\sigma_0(t) = \omega_{01}(t)\eta(t)$, is given by,

$$\lambda(t) = \sigma_0(t) - \frac{\dot{v}_I(t)}{2v_I(t)} + \omega_f(t) \begin{cases} \tanh(\varphi_f(t) - K) \\ \coth(\varphi_f(t) - K) \end{cases}$$
$$\omega_f(t) = \omega_{01}(t)v_I(t) \; ; \; \varphi_f(t) = \int \omega_f(t) dt$$
$$\lambda(t) = \sigma_0(t) - \frac{\dot{v}_{Im}(t)}{2v_{Im}(t)} + \omega_{fm}(t) \begin{cases} -\tan(\varphi_{fm}(t) - K) \\ \cot(\varphi_f m(t) - K) \end{cases} \quad (23)$$
$$\omega_{fm}(t) = \omega_{01}(t)v_{Im}(t) \; ; \; \varphi_{fm}(t) = \int \omega_{fm}(t) dt$$

The trigonometric forms assume that the intrinsic component is imaginary–namely, $jv_{Im}(t)$. These results provide powerful insight into the solution of LTV systems. Observe that integrating and exponentiating the eigenvalues produces the time domain solutions for the LTV systems to within a scale factor, A. For the case of real dynamic eigenvalues we have a time domain solution, y(t) = Ag(t)f(t), where,

$$g(t) = e^{\int (-\frac{\dot{v}_I(t)}{2v_I(t)} + \sigma_0(t))dt} = \frac{1}{\sqrt{v_I(t)}} e^{\int \sigma_0(t) dt}$$
$$f(t) = \begin{cases} \cosh(\varphi_f(t) - K) \\ \sinh(\varphi_f(t) - K) \end{cases} \quad (24)$$

For the case of complex dynamic eigenvalues we have a time domain solution where,

$$g(t) = \frac{1}{\sqrt{v_{Im}(t)}} e^{\int \sigma_0(t) dt} \; ; \; f(t) = \begin{cases} \cos(\varphi_{fm}(t) - K) \\ \sin(\varphi_{fm}(t) - K) \end{cases} \quad (25)$$

Note how in each case the intrinsic component of the RCE primitive pair essentially specifies the solution given a choice for K. By the appropriate choice of A and K, all possible boundary conditions in the time domain may be satisfied.

Note that the choice of A is a degree of freedom that may only be exploited in the time domain solution. However, K influences the choice of solution in the Riccati domain. Observe that the choice of K determines a specific boundary condition in the Riccati domain. Pursuing this point, we observe that the primitive pair of solutions are applicable only to a unique pair of initial conditions, which we may designate as $v_1(t_0) = v_R(t_0) + v_I(t_0)$ and $v_2(t_0) = v_R(t_0) - v_I(t_0)$. Using these conditions, the time domain solutions have f(t) equal to a purely exponential function. If a solution to the Riccati equation is desired with initial conditions not equal to either $v_1(t_0)$ or $v_2(t_0)$ then an appropriate member of the continuum of solutions must be used. In the case of real dynamic eigenvalues, given the fact that tanh() is a function with values only in the range (-1,1) for finite argument, and coth() is a function with values only in the range (- ∞, -1) and (1, ∞), all other initial conditions will be obtained only by one choice of tanh() or coth() with a unique choice of K.

## V. New Insights on Classic Theory

Before continuing the present discussion, some points regarding the RCE and its solutions can be made. To begin, if the intrinsic solution of the RCE becomes zero at any time, we have the paradoxical case where two independent solutions may propagate forward or reverse in time from a single boundary condition. Furthermore, if the intrinsic solution of the RCE becomes zero at an instant in time, say $t_0$, then $v_I(t)$ must have a factor of $(t - t_0)^n$, where n > 0. Given the relationship between the real and intrinsic components of the solution, then the real solution must include a term, given by $-n/[2\omega_{01}(t)(t - t_0)]$, plus others. Not only does this term become infinite, but $-n/2(t - t_0)$ goes to positive infinity as t increases to $t_0$. Hence, the integral of $\omega_{01}(t)v_I(t)$ will go to positive infinity at $t_0$, causing the associated time domain solution to be infinite. This scenario is impossible for the solution of a linear second order system with bounded parameters. Consequently, the intrinsic solution to the RCE, with bounded coefficients, may never become zero at finite time. In cases where the intrinsic solution goes to zero, with unbounded coefficients, the RCE solutions become equal at infinity, yielding a plausible intersection of independent solutions.

The restriction on a zero value for the intrinsic solution is valuable, effectively requiring boundedness of the parameters, $\omega_{01}(t)$ and $\omega_{02}(t)$. However, a zero value for $\omega_{01}(t)$ is also problematic. Many systems of interest have an RCE with a value of 1 for $\omega_{01}(t)$. For example, any system described by the scalar differential equation in (3), exemplified by the engineering examples given below in section V, have $\omega_{01}(t) = 1$. Nevertheless, one may certainly find cases where $\omega_{01}(t)$ changes sign one or more times for times of interest. While this may introduce finite escape events for $v_R$, it may not create, on its own, infinite values for the dynamic eigenvalues, since they are scaled by $\omega_{01}(t)$. However, the product $\omega_{01}(t)v_I(t)$, given $v_I(t)$ never equals zero, will change sign and its integral, $\varphi_f(t)$, will not be monotonic, causing potential finite escape events not correlated simply to the primitive solutions. While this scenario cannot be ruled out for potential systems of interest, such systems seem unusual and will be considered contrived in the present context.

Another observation regarding the RCE parameters is worth noting. Namely, when $\omega_{01}(t)$ and $\omega_{02}(t)$ are nonzero and have opposite sign then the RCE is inherently unstable and predisposed to exhibit a finite escape time solution. This is because when $-\omega_{01}(t)$ and $\omega_{02}(t)$ have the same sign, the derivative of the solution, $v(t)$, is globally positive or negative. If these parameters have the same sign, the RCE can exhibit solutions that are asymptotically stable.

**V.a RCE Phase Plane**

Let us consider the case where $\omega_{01}(t)v_I(t)$ is positive real for all time, and $\omega_{01}(t)$ and $\omega_{02}(t)$ are bounded. Given the above comments, this constraint will apply to well behaved (bounded and not contrived) systems, given a purely real primitive solution to the RCE. Note that in designating $v_R(t) \pm v_I(t)$ as primitive solutions there is no loss of generality in this positive assumption. Then $\varphi_f(t)$, the integral of $\omega_{01}(t)v_I(t)$, grows positively without bound and all solutions will converge toward $v_R(t) + v_I(t)$ as $t \to \infty$, since both tanh() and coth() approach 1. However, if initial conditions dictate that the choice of K makes the argument of the hyperbolic functions negative at the starting time of interest, then only the tanh() function will remain bounded as time goes forward, while the coth() functional form guarantees finite escape time. This is an interesting observation in that the primitive pair of solutions are exactly the attractor (fixed point in the limit) and the separatrix in a functional sense for the given Riccati equation.

This fundamental relation between the separatrix and the attractor has not been identified before. A simple example illustrates this idea. Suppose the Riccati equation is $\dot{\nu} = -\nu^2 + 4$. The primitive solution pair is $v = \pm 2$. It is easily shown that $v = -2$ is the separatrix, since any initial condition less than -2 produces the solution 2coth(2t - K), with a positive value of K, guaranteeing a finite escape time solution at t = K/2. Any initial condition above -2 spawns either 2tanh(2t - K) or 2coth(2t - K), with a negative value of K, as a solution. Both of these solutions converge towards the attractor, $v = 2$, whose limit in time is 2, the fixed point. Note that even solutions with a finite escape time event will converge to the attractor following the event. A time varying Riccati equation, having a primitive solution where $\omega_{01}(t)v_I(t)$ is positive real for all time obeys the same behavior, with time varying separatrix and attractor.

An interesting aspect of the RCE solution converging to the attractor is that one may view the attractor solution as a steady state response, where initially there is what may be called an RCE transient, as the RCE solution converges to the attractor. This phenomenon is not consistent with the usual interpretation of the solution to a linear homogeneous differential equation, whose response is thought of as the transient of a complete response to an input. However, the RCE transient concept may well be valuable in LTV system analysis. To appreciate its connection even to LTI systems, consider the above example where the primitive solution pair is $v = \pm 2$. This LTI system example would typically produce a pair of exponentials, $e^{2t}$ and $e^{-2t}$, in the time domain. A given set of initial conditions would excite both exponential terms; however, the growing exponential would rapidly swamp the decaying exponential, leaving what would appear to be the single exponential "steady state" result. The time during which both exponentials were noticeably present could be referred to as the "RCE transient" period. This interval may be more valuable to understand in LTV engineering applications.

An alternate argument holds for the case of complex primitive solutions to the RCE, where we assume well behaved (bounded and not contrived) systems. In that case, real boundary conditions may never result in a primitive solution. Only the tan() or cot() forms, which are intrinsically equivalent, will apply to satisfy arbitrary real initial conditions. Furthermore, all of these solutions will exhibit repeated finite escape time, corresponding to repeated zero crossings in the time domain solutions. As implied above, an allowable finite escape time event must approach negative infinity as time increases. That way its integral approaches negative infinity, so that the time domain solution will approach zero. Notice that -tan and cot functions then flip to positive infinity with negative derivative, which will bring the time domain solution away from zero. This behavior is completely predicted by the solutions of (25). Recall that in LTI systems complex conjugate eigenvalues result in sinusoidal time domain solutions. In this case the phase plane generally has neither a separatrix nor an attractor that are global. Exceptions to this behavior may occur in cases where dynamic eigenvalues transition between positive and negative values or real and complex values at finite time. Other interesting behavior may occur in a RCE corresponding to a periodic system. This is discussed later.

One final note is that the presence of an attractor for a system may be detected in a variety of ways, both analytically and in simulation, based on the methodology presented here. The presence of an attractor guarantees that the primitive solution is real. This knowledge leads to an accurate guess in finding the primitive solution as described in section III.

**V.b Polynomial Solutions**

It is noted that the general form of solutions for the Riccati equation given here is not obviously consistent with polynomial solutions given by example and in the literature [e.g.8] for explicitly solvable Riccati equations with boundary conditions imposed. For example, the Riccati equation given by, $\dot{\nu} = -\nu^2 + 2/t^2$, has 2 known simple solutions, z = -1/t and 2/t. These two solutions may satisfy only two possible boundary conditions at a given time. The general solution is $(2t^3 - C)/t(t^3 + C)$, where C may be chosen to satisfy any specific initial condition at a given time. While this appears to be a different form than given above, in fact it is completely consistent with the results here. Specifically, given the reduced form of this Riccati equation, the known solutions, z = v = -1/t and 2/t, lead to the conclusion that $\nu_R = 1/2t$ and $\nu_I = 3/2t$, defining the primitive solution, $\nu_p(t)$. Therefore, we may find the general solution, $\nu(t)$.

$$\nu_p(t) = \nu_R(t) \pm \nu_I(t) \ ; \ \nu_R(t) = \frac{1}{2t} \ , \ \nu_I(t) = \frac{3}{2t}$$
$$\nu(t) = \frac{1}{2t} + \frac{3}{2t} \begin{cases} \tanh(\int \frac{3}{2t} - K) \\ \coth(\int \frac{3}{2t} - K) \end{cases} \quad (26)$$

The hyperbolic functions may be transformed to polynomial ratios due to the form of $\nu_I(t)$, with the help of the results in (15).

$$e^{(\int \frac{3}{2t} - K)} = \frac{t^{3/2}}{\sqrt{C}} \ ; \ C = e^{2K}$$
$$\Rightarrow \tanh(\int \frac{3}{2t} - K) = \frac{t^3 - C}{t^3 + C} \quad (27)$$
$$\nu(t) = \frac{1}{2t} + \frac{3}{2t} \frac{t^3 - C}{t^3 + C} = \frac{2t^3 - C}{t(t^3 + C)}$$

Notice that this solution, corresponding to the tanh() option, assumes that C must be positive due to its relation to K. However, if the coth() option is used, the solution above applies with assumed negative values for C. Hence, the solution for this problem is a special case of the general solution given here. To the author's knowledge, no known polynomial solution to the differential Riccati equation fails to fit the solution form of this work.

**V.c The Periodic case:**

Consider the special case where a second order LTV system has a periodic state matrix–that is, A(t) = A(t + T). This can only be the case if each element of A(t) is periodic. It has been shown [1] that a fundamental matrix for any second order LTV system is given by,

$$\phi(t) = V(t) \begin{bmatrix} e^{\int \lambda_1(t)dt} & 0 \\ 0 & e^{\int \lambda_2(t)dt} \end{bmatrix}$$
$$V(t) = \begin{bmatrix} 1 & 1 \\ \nu_1(t) - \alpha(t) & \nu_2(t) - \alpha(t) \end{bmatrix} \quad (28)$$

where $\nu_1(t)$ and $\nu_2(t)$ are the solutions to the Riccati characteristic equation associated with the system. Given the periodicity of A(t), then $\omega_{01}(t)$, $\omega_{02}(t)$, $\alpha(t)$, and $\sigma_0(t)$ are also periodic. Hence, the eigenvector matrix, V(t), and the dynamic eigenvalues, $\lambda_1(t)$ and $\lambda_2(t)$, are periodic if and only if the solutions of the Riccati characteristic equation are periodic.

Given periodic dynamic eigenvalues, their respective integrals must be the sum of a zero mean periodic function and a constant times time, where the constant, $r_i$, is the respective constant, which may be called the DC value, associated with each of the dynamic eigenvalues. Thus, we may write the fundamental matrix of the LTV system as follows.

$$\phi(t) = V(t) \begin{bmatrix} e^{p_1(t)} & 0 \\ 0 & e^{p_2(t)} \end{bmatrix} \begin{bmatrix} e^{r_1 t} & 0 \\ 0 & e^{r_2 t} \end{bmatrix} \quad (29)$$
$$\int \lambda_{1,2}(t) dt = p_{1,2}(t) + r_{1,2} t \ ; \ r_{1,2} = \frac{1}{T} \int_0^T \lambda_{1,2}(t) dt$$

Continuing, we have,

$$\phi(t) = Q(t) W \begin{bmatrix} e^{r_1 t} & 0 \\ 0 & e^{r_2 t} \end{bmatrix} W^{-1} = Q(t) e^{Rt}$$
$$Q(t) = V(t) \begin{bmatrix} e^{p_1(t)} & 0 \\ 0 & e^{p_2(t)} \end{bmatrix} W^{-1} \ ; \ R = W \begin{bmatrix} r_1 & 0 \\ 0 & r_2 \end{bmatrix} W^{-1} \quad (30)$$

W is any non-singular square constant matrix. Clearly, Q(t) is periodic and $e^{RT}$ is a matrix, typically called the monodromy matrix [9], with Floquet exponents given by $r_1 T$ and $r_2 T$. If the solutions to the RCE are not periodic, then Q(t) cannot be periodic due to the given form of V(t) and the fact that the exponentials of $p_i(t)$ will both be aperiodic. However, Floquet's theorem guarantees that Q(t) must be periodic. Hence, the solutions to the RCE must be periodic any time $\omega_{01}(t)$ and $\omega_{02}(t)$ are periodic.

This has never before been proven [e.g.,10,11]. However, it should be noted that Floquet's theorem assumes a steady state solution. Given a periodic system, arbitrary initial conditions will not yield a periodic solution, at least not initially. In the case of real primitive solutions for the RCE, any solution will eventually equal the attractor which must obey Floquet's theorem. In the short term, the hyperbolic functions may preclude periodicity. To be specific, the intrinsic solution must be periodic with a DC average. The relation between the real and intrinsic solutions then guarantees that the real solution will be periodic with zero DC average. A nonzero DC average guarantees that the hyperbolic functions will eventually go to ±1, yielding the eventual primitive solution. The DC component will set the Floquet exponents [9] to be ±T times the DC average of the intrinsic primitive solution. It is noted that for second order systems, then, the product of the characteristic multipliers must be one if η(t) equals zero, as in the case of the Matthieu equation for example.

Given an imaginary intrinsic primitive solution to the RCE, comments regarding the DC average are unchanged except that the DC average is an imaginary number. The real component still has a zero DC average. The intrinsic solutions to the RCE are now complex conjugate functions. The observable real solutions to the RCE will involve trigonometric functions that will have periodic arguments plus a linear function of time. Hence, we expect likely aperiodic repetitive finite escape time solutions, unless the DC average of the intrinsic solution has an appropriate value. Rather complex time domain behavior can result from complex periodic RCE solutions. This will be discussed further below in relation to the well studied Mathieu equation.

**V.d Engineering Examples**

In this section, three well studied systems are considered in the context of the present analysis. In each case, aspects of the systems or their current analysis appear to be at odds with the analysis here. But in each case, it is shown that the present analysis is validated, and in fact predicts and verifies the well known results.

Let us begin with the case of Bessel functions. The classic Bessel differential equation, parametrized with integers, admits solutions that are linear combinations of so-called Bessel functions of the first and second kind [12]. For small values of time, these two types of time domain solutions are drastically different; however, for large time they are extremely similar. To be specific, the Bessel differential equation considered here is

$$t^2 \ddot{y}(t) + t\dot{y}(t) + (t^2 - N)y(t) = 0$$
$$\Rightarrow \ddot{y}(t) + \frac{1}{t}\dot{y}(t) + (1 - \frac{N^2}{t^2})y(t) = 0 \quad (31)$$

We will allow the singularity created by dividing by $t^2$, as is normal in practice. We have the state equations as in (4) given by,

$$\begin{bmatrix} \dot{x}_1(t) \\ \dot{x}_2(t) \end{bmatrix} = \begin{bmatrix} 0 & 1 \\ \frac{N^2}{t^2} - 1 & -\frac{1}{t} \end{bmatrix} \begin{bmatrix} x_1(t) \\ x_2(t) \end{bmatrix} \quad (32)$$

Note that this is a case where $\omega_{01}(t) = 1$. The value of $\omega_{02}(t)$, however, becomes infinite at t = 0, allowing an infinite time domain solution. Using the present notation, we may find the related RCE and dynamic eigenvalues using the procedures in [1] restated earlier in (5) and (6), given by,

$$\dot{\nu}(t) = -\nu^2(t) + \omega_{02}(t) \; ; \; \omega_{01}(t) = 1$$
$$\omega_{02}(t) = \frac{4N^2 - 1}{4t^2} - 1; \; \lambda(t) = \nu(t) - \frac{1}{2t} \quad (33)$$

Note that the parameter, N, is assumed here to be a positive integer which generates the families of solutions associated with Bessel functions of the first and second kind. It is interesting to note that asymptotically it would appear that the eigenvalues of this equation should be real when t → 0 and complex as t → ∞. Nevertheless, if one back propagates the solution starting from a good guess at the complex solution for large t, the actual solution is complex for all t > 0. Using 5 as the value of N, as shown in Fig 1, the complex primitive solution rapidly approaches a real function for small time, but the intrinsic solution remains positive imaginary for all t > 0. Consequently, for all positive time the solution to the Bessel RCE possesses a purely complex solution. As $\nu_I(t)$ approaches zero at t = 0, the real component approaches infinity, as expected, avoiding the situation where a pair of primitive solutions unify at a finite value. The conjugate of this solution also solves (using forward propogation in simulation) the RCE, verifying the existence of the complementary pair.

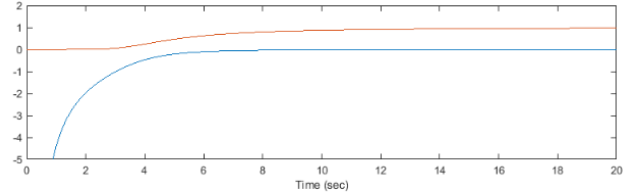

Figure 1: Real (blue) and imaginary (red) parts of the primitive RCE solution.

The primitive solutions are unobservable given any purely real boundary condition. Instead, we find a pair of observable real solutions using the methodology of this work. Because the intrinsic solution is imaginary we use the results in (18) to express the solution. A value of K = π/2 in these expressions allows classic initial conditions for Bessel functions of the first and second kind to be satisfied for the complementary pair of solutions. This choice for K simply interchanges the -tan and cot functions for the complementary solutions of (18). Using the back propagated primitive solution shown in Fig 1 as $\nu_{Im}(t)$, and whose integral is $\varphi_{fm}(t)$, both analytically based solutions for ν(t) from (18) are plotted in Fig 2. In addition, the actual simulated solution for the RCE, with initial condition taken from the selected first solution, $\nu_R(t) + \nu_{Im}(t)\cot(\varphi_{fm}(t))$ at t = 0.01 sec, is shown to exactly overlap the analytic one up to the first finite escape time and beyond. The analytic and simulated solutions clearly show the repetitive finite escape time solutions for larger times.

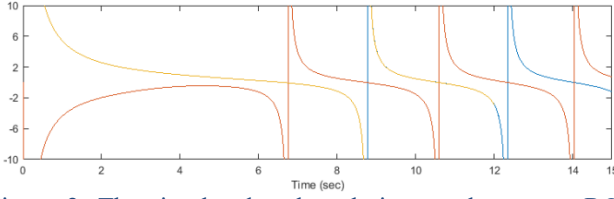

Figure 2: The simulated real analytic complementary RCE solutions (blue and red) and the actual simulated solution (yellow) overlapping the blue solution.

The integral of the pair of analytic solutions, offset by 1/2t as indicated in (33) reveals the respective associated solutions in the time domain. The calculations are as follows:

$$\begin{aligned}
\lambda_1(t) &= -\frac{1}{2t} - \frac{\dot{\nu}_{\text{Im}}(t)}{2\nu_{\text{Im}}(t)} + \omega_{fm}(t)\cot(\varphi_{fm}(t)) \\
\lambda_2(t) &= -\frac{1}{2t} - \frac{\dot{\nu}_{\text{Im}}(t)}{2\nu_{\text{Im}}(t)} - \omega_{fm}(t)\tan(\varphi_{fm}(t)) \\
y_1(t) &= e^{\int \lambda_1(t)dt} = \frac{A_1}{\sqrt{t\nu_{\text{Im}}}}\sin(\varphi_{fm}(t)) \\
y_2(t) &= e^{\int \lambda_2(t)dt} = \frac{A_1}{\sqrt{t\nu_{\text{Im}}}}\cos(\varphi_{fm}(t))
\end{aligned} \quad (34)$$

When the solutions, $y_1(t)$ and $y_2(t)$, are scaled by 0.8 and -0.8, respectively, they match the Bessel functions of the first and second kind computed by MATLAB using besselj and bessely functions. Clearly the theory presented here explains the solutions, and intuitively identifies Bessel functions of the first and second kind as being associated with complementary solutions of the RCE for the Bessel differential equation. The complementary forms for these solutions in (34) are also innovative compared to classic work.

These findings are not coincidental. The classic well-studied quantum harmonic oscillator [13,14] is another example that exemplifies the validity of the approach offered here. The core differential equation associated with that system and its RCE are given below, where N is an odd positive integer.

$$\begin{aligned}
\ddot{y}(t) + (N - t^2)y(t) &= 0 \; ; \; \dot{\nu}(t) = -\nu^2(t) + \omega_{02}(t) \\
\omega_{01}(t) &= 1 \; ; \; \omega_{02}(t) = t^2 - N
\end{aligned} \quad (35)$$

Analytic solutions have been found for this system using Hermite polynomials. Analogously to the Bessel case just presented, one may guess that the dynamic eigenvalues are complex for small time and real for large time. However, by making a good guess at small time that is complex, simulation confirms that the solution to the RCE is complex for all finite time greater than or equal to zero. And as in the last case, using the imaginary intrinsic component one may use the present theory to exactly specify the pair of observable real solutions that spawn the well recognized solutions to this problem. Using the simulated primitive solution, with calculations similar to those above, the solutions predicted by the present analysis using the tan option exactly (to within computational accuracy) overlap the well established known solutions for the wavefunctions of bound eigenstates.

Another noteworthy case is the well studied Matthieu equation [15] whose time domain equation and RCE are given here.

$$\begin{aligned}
\ddot{y}(t) - (a_0 + q\cos(t))y(t) &= 0 \\
\dot{\nu}(t) &= -\nu^2(t) + a_0 + q\cos(t)
\end{aligned} \quad (36)$$

For all values of the parameters, $a_0$ and q, the RCE must possess periodic solutions as per the above discussion. This is despite the fact that the time domain solutions exhibit a broad spectrum of behavior, which is well known [9]. Note that $\eta(t) = \alpha(t) = 0$ and $\omega_{01}(t) = 1$ for this system; therefore, the dynamic eigenvalues are simply given by a complementary pair of solutions to the RCE. Given this and the earlier discussion, the dynamic eigenvalues corresponding to steady state solutions are given by a periodic zero average (DC value) real component, $\nu_R(t)$, plus or minus a periodic intrinsic component, $\nu_I(t)$. The intrinsic component may have a nonzero DC value. If $\nu_I(t)$ is real the DC value will cause a pair of time domain solutions where one grows and the other decays in time, as the RCE solution converges toward an attractor. Following the RCE transient, the "steady state" will be an exponentially growing function in the time domain, indicative of an unstable system solution. If $\nu_I(t)$ is imaginary the DC value will cause an apparently aperiodic mix of sinusoids in the time domain in any short window of time, unless the DC value takes on a special value. However, these time domain solutions will remain bounded whether periodic or not. Note that these observations clearly imply that the Floquet exponents found using (30) will either be a pair of real constants or imaginary constants. In the real case, time domain solutions will be unstable and the Riccati domain solution will tend toward an attractor. In the imaginary case, time domain solutions will be bounded (stable) and often aperiodic.

It is noteworthy that these observations are at odds with some aspects of well accepted Floquet theory. In particular, when a Mathieu system time domain solution is an exponentially weighted periodic signal of period 2T, accepted current analysis [9,15] exploits the fact that the periodic component, $y_p(t)$, is anti-symmetric with $y_p(t + T) = -y_p(t)$. In this case, the otherwise real Floquet exponent is augmented with ±jπ/T endowing the monodromy matrix with the anti-symmetric property, and preserving the mathematical accuracy of the Floquet theorem. This is inconsistent with the theory proposed here.

In order to explore this and again test the work here, a series of simulations were run with q equal to 1 and $a_0$ varied over a range of values from positive to negative. In each case the initial conditions were varied in order to see the evolution of solutions and any eventual steady state behavior. The inspection of many solutions with carefully chosen initial conditions allowed quite accurate estimates of real and intrinsic components, including any transient behavior. It was observed that for positive values of $a_0$, greater than approximately 0.4, all cases showed periodic, T = 2π, RCE solutions with a clear attractor and separatrix. In each case there was a DC average that could be directly associated with a real intrinsic component that remained greater than zero for all time. In the time domain the associated solutions were exponentially weighted periodic, T = 2π, functions. The exponential weighting was completely consistent with the DC average of the intrinsic component. For negative values of $a_0$, less than approximately -1.5, no attractor was observed and solutions to the RCE showed repeated finite escape time, indicative of a complex RCE solution, whose intrinsic component had an imaginary DC average. For example, with $a_0$ = -3, careful choice of a complex initial value for the RCE produced a periodic, T = 2π, complex solution whose intrinsic component was imaginary with a DC average of 1.72, which was very close to $\sqrt{3}$ as expected. This was consistent with the associated time domain solutions that appeared as a bounded combination of sinewaves. Careful choice of $a_0$, using values taken from the literature, showed a bounded harmonic combination of sinewaves that was periodic with T = 2π.

For cases where $a_0$ was between approximately -1.5 and 0.4, the results were more complicated, but still consistent with the current analysis. For approximately -0.55 < $a_0$ < 0.3, a clear attractor was observed, indicating a real primitive solution pair; however, the intrinsic solution showed periodic finite escape time events, with periodic finite escape time events for the real component as well. Only the real component crossed zero, and there was a clear DC average to the intrinsic component. The composite primitive solution to the RCE was periodic, T = 2π, and included a single finite escape time event per period. This was crucial, since the steady state time domain solution was an exponentially weighted periodic function with a period of 2T. The finite escape time of the RCE solution explained one zero crossing per T seconds, while the periodicity of the RCE solution explained that the periodic function, $y_p(t)$, had anti-symmetric half cycles, such that $y_p(t + T) = -y_p(t)$, for all t following the RCE transient. The exponential weighting was consistent with the estimated DC average of the intrinsic component. Solutions with $a_0$ in the range between approximately -1.5 and -1 showed an attractor, and with $a_0$ in the range between approximately -1 and -0.6 showed no attractor. All solutions with attractors showed unbounded periodic behavior in the time domain. Those without attractors showed a bounded generally aperiodic mix of sinusoids in the time domain, indicative of an RCE imaginary intrinsic solution equal to a constant (DC average) plus a periodic function. In the time domain these solutions were often aperiodic on the time scale of simulation, likely due to a DC average not simply related to 2π/T.

While this examination of the Matthieu equation was inadequate to make a concrete mathematical judgment, it certainly tended to confirm all the expectations stated above based on the theory proposed here. Regarding cases with 2T periodic solutions, it is acknowledged that adding ±jπ/T to the Floquet exponents correctly allows the statement of Floquet's theorem to be accurate. But it artificially multiplies the fundamental matrix, which is an anti-symmetric function with a period of 2T, by another anti-symmetric function with a period of 2T so that the product of the functions has a period of T. The modulation of the fundamental matrix is undone using the monodromy matrix whose Floquet exponents are augmented with ±jπ/T. The net result is to hide the fact that all truly complex eigenvalues for the monodromy matrix correspond to complex dynamic eigenvalues which, in the case of the Matthieu equation, guarantee stable bounded solutions in the time domain, at least according to the theory here. This merits further study.

**VI. Discussion**

This paper has offered a systematic look at the solutions to the Riccati differential equation, referred to here as the Riccati characteristic equation (RCE), owing to its intimate relation to all possible linear time-varying second order systems, as discussed in the companion paper [1]. For the first time all solutions to this equation have been shown to be given by a general standard form generated by an intrinsic solution that completely specifies a primitive pair of solutions that may be found for all Riccati equations. All Riccati equations connected to practical systems possess real solutions, despite the fact that a subset possesses complex primitive solutions. Due to its relation to LTV systems the solutions to the RCE determine pairs of dynamic eigenvalues, which when integrated and exponentiated provide the solutions "in the time domain" to the associated LTV systems. As a consequence of this, finite escape time solutions to the Riccati equation are not disqualifying, but rather imperative, since each associated time domain solution will likely exhibit zero crossings that require the RCE solutions to go to negative infinity with increasing time. The general form of solutions offered in this work gracefully allows repetitive finite escape time solutions. This general form also gracefully predicts an attractor (with a limit point) and separatrix in one class of cases, and parallel flow with finite escapes in another.

Perhaps the most important aspect of this work is the clear unification of all known solutions, from polynomial to Bessel and periodic functions. The paradigm of thinking of complementary pairs of solutions as shown here allows for a generalization to the study of linear systems that has never been available. Consequently, the perspective offered by the work here provides completely new insight into the behavior of well studied engineering systems and the Riccati differential equation itself. For example, Floquet's theorem has been shown to prove that a pair of periodic solutions to the Riccati differential equation exist whenever the equation has periodic parameters, and that the commonly accepted Floquet theory analysis may benefit from a modification presented here for the first time. It is hoped that the insight provided by the work here may improve the understanding and teaching of the Riccati differential equation and general linear time varying systems.